\documentclass[11pt]{article}
\usepackage{amssymb}
\usepackage{amsmath}

\newtheorem{theo}{Theorem}
\newtheorem{prop}{Proposition}
\newtheorem{lemm}{Lemma}
\newtheorem{coro}{Corollary}
\newtheorem{rema}{Remark}
\newtheorem{Defi}{Definition}
\newtheorem{ex}{Example}
\newtheorem{conj}{Conjecture}

\newcommand{\cqfd}
{%
\mbox{}%
\nolinebreak%
\hfill%
\rule{2mm}{2mm}%
\medbreak%
\par%
}
\newfont{\gothic}{eufb10}
\date{\empty}
\begin{document}
\title{Intrinsic pseudovolume forms and $K$-correspondences}
\author{Claire Voisin\\ Institut de math{\'e}matiques de Jussieu, CNRS,UMR
7586}

\maketitle \setcounter{section}{-1}
\section{Introduction} In recent years, the notion of
$K$-equivalence has appeared in several contexts, like motivic
integration \cite{denef}, McKay correspondence \cite{batyrevdais}
and derived category of coherent sheaves on varieties
\cite{kawamata}, \cite{wang}. A $K$-equivalence between two
algebraic varieties $X$ and $Y$ is a birational map
$\phi:X\dashrightarrow Y$ whose graph $\Gamma_\phi\subset X\times
Y$ admits a desingularization
$$\tau:Z\rightarrow \Gamma_\phi$$
such that $f^*K_X$ and $ g^*K_Y$ are linearly equivalent, where
$f=pr_1\circ \tau,\,g=pr_2\circ \tau$. Equivalently, the two
ramification divisors  should satisfy : \begin{eqnarray}
\label{eqnabove}R_f=R_g. \end{eqnarray}In this paper, we start the
study of what we call $K$-(iso)correspondences  between smooth
varieties or complex manifolds $X,\,Y$ of the same dimension,
which are  graphs of multivalued maps, or analytic subsets in the
product $X\times Y$, generically finite over each factor, such
that any desingularization $\widetilde \Sigma$ satisfies with the
notations above
 the condition (\ref{eqnabove}) or, in case of
$K$-correspondences, the weakened condition
$$R_f\leq R_g.$$
Hence we simply forget the condition that the degree of the graph
over $X$ and $Y$ should be $1$. Such $K$-isocorrespondences appear
naturally in the McKay situation (cf section 2).

Our main result proved in section 2, is the fact that many
$K$-trivial projective varieties carry a lot of
self-$K$-isocorrespondences $\Sigma\subset X\times Y$ satisfying
the condition that $deg\,{pr_1}_{\mid \Sigma}\not=
deg\,{pr_2}_{\mid \Sigma}$. With the notations above, one can see
easily that this last condition is equivalent to the equality of
volume forms on $\widetilde \Sigma$
$$ f^*\Omega_X=\lambda g^*\Omega_X,$$
where $\lambda\not=1$ is a real number, and $\Omega_X$ is the
canonical volume form of $X$. Equivalently
$$f^*\omega_X=\mu g^*\omega_X$$
where $\omega_X$ is any generator of $H^0(X,K_X)$, and $\mu$ is a
complex number of modulus $\not=1$. Because of this property,
these self-$K$-isocorrespondences look like the multiplication by
an integer in an abelian variety.

Section 3 discusses potential applications of this result to the
study of intrinsic pseudovolume forms on complex manifolds (see
\cite{kobayashi}). Kobayashi  and Eisenman have introduced an
intrinsic pseudovolume form $\Psi_X$ on any complex manifold $X$,
which is computed using all holomorphic maps from a polydisk
$D^n,\,n=dim\,X$ to $X$. Here we introduce modified intrinsic
pseudovolume forms
$$\Phi_{X,an}\leq \Phi_X\leq \Psi_X$$
which are defined essentially by replacing holomorphic maps with
holomorphic $K$-correspondences in the Eisenman-Kobayashi
definition. We show that on one hand, the following theorem, due
to Green-Griffiths \cite{grigre} and Kobayashi-Ochiai
\cite{kobaochiai} , still holds for the pseudovolume form
$\phi_{X,an}$ :
\begin{theo} If $X$ is a projective variety which is of general
type, $\Phi_{X,an}>0$ on a dense Zariski open set of $X$.
\end{theo}

On the other hand we show that $\Phi_X$  is equal to $0$ for a
large number of families of $K$-trivial varieties listed in
section 3, and also for varieties which are fibered with fiber of
these types. This gives us a weak version of the Kobayashi
conjecture (i.e. the converse to the
Green-Griffiths-Kobayashi-Ochiai theorem) for the modified
pseudovolume form $\Phi_X$.

Section 4 is devoted to a few supplementary results, remarks and
questions concerning $K$-isocorrespondences. In particular, we
provide (cf corollary \ref{coro}) new examples of $K$-trivial
varieties satisfying Kobayashi's conjecture \ref{conjko}.

{\bf Acknowledgements}. I wish to thank Fr{\'e}d{\'e}ric Campana, for
communicating to me his preprint \cite{campana}, and for asking
interesting questions which on related topics, which led me to
work on this subject.
\section{$K$-correspondences}
In this section, we introduce and discuss the notions of
$K$-correspondences and $K$-isocorrespondences, which are
straightforward generalizations of the so-called $K$-ordering and
 $K$-equivalence in birational geometry (cf \cite{kawamata}, \cite{wang}).
 We assume that $X$ and $Y$ are smooth complex manifolds of
 dimension $n$.
 \begin{Defi} A $K$-correspondence from $X$ to $Y$ is a reduced
 $n$-dimensional closed analytic subset $\Sigma\subset X\times Y$,
 such
 that on each irreducible component of $\Sigma$, the projections
 to $X$ and $Y$ are generically of maximal rank, and satisfying the
 following two conditions :
 \begin{enumerate}
 \item The restriction ${pr_1}_{\mid \Sigma}$ is proper.
 \item \label{itemb}Let $\widetilde\Sigma\stackrel{\tau}{\rightarrow}\Sigma$ be
 a desingularization, and let
 $$f:=pr_1\circ\tau:\widetilde\Sigma\rightarrow X,\,g=pr_2\circ
 \tau:\widetilde\Sigma\rightarrow Y.$$
 Then we have the inequality of ramification divisors on
 $\widetilde\Sigma$ :
 $$R_f\leq R_g.$$
 \end{enumerate}

 \end{Defi}
 Note that  property \ref{itemb} has to be checked on
 one desingularization, and then will be satisfied by all
 desingularizations, as a standard argument shows.

 A holomorphic map $\phi$ from $X$ to $Y$ leads to a correspondence,
 obtained by taking
 the graph of $\phi$. It turns out that
 $K$-correspondences behave with many respects as maps. Their main
 common feature with ordinary maps is the fact that for any
 desingularization $\tau:\widetilde\Sigma\rightarrow \Sigma$ as
 above, we get a  natural inclusion
 $$g^*K_Y\subset f^*K_X,$$
 as subsheaves of $K_{\widetilde\Sigma}$. We also have the
 following important fact :
 \begin{prop} $K$-correspondences can be composed. More precisely,
 if
 $\Sigma\subset X\times Y$ and $\Sigma'\subset Y\times Z$ are
 $K$-correspondences, then define
 $\Sigma'\circ\Sigma$ to be the union of the components
 of $p_{13}(p_{12}^{-1}(\Sigma)\cap p_{23}^{-1}(\Sigma'))$ on
 which the projections to $X$ and $Z$ are generically of maximal
 rank. Then $\Sigma'\circ\Sigma$ is a $K$-correspondence.
 \end{prop}
 {\bf Proof.} Note first that the properness of the first
 projections on
 $\Sigma$ and $\Sigma'$ implies that
 $p_{13}(p_{12}^{-1}(\Sigma)\cap p_{23}^{-1}(\Sigma'))$
 is a closed analytic subset of
 $X\times Z$. It also shows that the projection to $X$ is proper
 on this analytic subset. Finally it is easy to see that a
 component which is generically of maximal rank over both
 $X$ and $Z$ must be of dimension $n$.

 Next let $\widetilde\Sigma,\,\widetilde\Sigma'$ be
 desingularizations of $\Sigma,\,\Sigma'$. Denote by $f,\,g$ the
 maps from $\widetilde\Sigma$ to $X$ and $Y$, and by
 $f'$ and $g'$ the maps from $\widetilde\Sigma'$ to $Y$ and $Z$.
 Let $\Sigma''$ be a component of
 $\widetilde\Sigma\times_Y\widetilde\Sigma'$ on which
 the maps $F:=f\circ \phi$ and $G:=g'\circ\psi$ are generically of
 maximal rank. Here $\phi:\Sigma''\rightarrow\widetilde\Sigma$ and
 $\psi:\Sigma''\rightarrow\widetilde\Sigma'$ are the two natural
 maps. Choose a desingularization $\widetilde\Sigma''$ of
$ \Sigma''$. Let now $\sigma\in \widetilde\Sigma''$ and let
 $$x=F(\sigma),\,z=G(\sigma),\,y=g\circ\phi(\sigma)=f'\circ\psi(\sigma).$$
 Let $\omega_x$ be a holomorphic $n$-form which generates $K_X$
 near $x$, and similarly choose $\omega_y$ near $y$ and $\omega_z$ near
 $z$. Then property \ref{itemb} says that we have the following
 equality of
 $n$-forms on $\widetilde\Sigma$ and $\widetilde\Sigma'$ respectively:
 $$g^*\omega_y=\chi\cdot f^*\omega_x,\,
 {g'}^*\omega_z=\chi'\cdot {f'}^*\omega_y,$$
 where $\chi$ is a holomorphic function on
 $\widetilde\Sigma$ and $\chi'$ is a holomorphic function on
$ \widetilde\Sigma'$, defined respectively on the inverse image in
$\widetilde\Sigma$ of a neighbourhood of $(x,y)$ in $X\times Y$
and on the inverse image in $\widetilde\Sigma'$ of a neighbourhood
of $(y,z)$ in $Y\times Z$.

 Pulling-back, via $\phi,\,\psi$ respectively, these equalities to
$\widetilde\Sigma''$ now  gives :
$$\chi\circ\phi \cdot F^*\omega_x=\phi^*(g^*\omega_y),$$
$$G^*\omega_z=\chi'\circ\psi\cdot\psi^*({f'}^*\omega_y).$$
Then, using $g\circ\phi=f'\circ\psi$, we conclude that
$\phi^*(g^*\omega_y)=\psi^*({f'}^*\omega_y)$ and hence
$$G^*\omega_z=\chi'\circ\psi\cdot\chi\circ\phi\cdot F^*\omega_x$$
as $n$-forms on $\widetilde\Sigma''$, where $\chi'\circ\psi\cdot
\chi\circ\phi$ is a holomorphic function on $\widetilde\Sigma''$.
 \cqfd
Our next definition is the following :
\begin{Defi} A $K$-isocorrespondence between $X$ and $Y$ is
a $K$-correspondence $\Sigma$ from $X$ to $Y$ such that $^t\Sigma$
is a $K$-correspondence from $Y$ to $X$, where $^t$ means the
image under the natural isomorphism $X\times Y\cong Y\times X$.
\end{Defi}
In other words, $\Sigma$ has to satisfy the properties that the
two projections $pr_i$ on $X$ and $Y$ are proper on $\Sigma$ and
that if $\tau:\widetilde\Sigma\rightarrow\Sigma$ is a
desingularization, with $f=pr_1\circ\tau,\,g=pr_2\circ \tau$, we
have now the equality
$$R_f=R_g.$$
$K$-isocorrespondences look like isomorphisms with certain
respects. The most important point for us will be the fact that
with the notations above, a $K$-isocorrespon-\\dence induces a
canonical isomorphism
$$f^*K_X\cong g^*K_Y.$$
Indeed they are equal as subsheaves of $K_{\widetilde\Sigma}$.

With the same arguments as before, one shows that
$K$-correspondences can be composed.
\begin{ex} If $f:X\rightarrow Y$ is a proper {\'e}tale map,
the graph of $f$ and its transpose are $K$-isocorrespondences.
\end{ex}
\begin{ex} If $G$ is a finite group acting on X, in such a way
that the stabilizer $G_x$ acts via $SL(n)$ on the tangent bundle
$T_{X,x}$ at each point $x$ of $X$, the quotient $X/G$ has
Gorenstein singularities. If $x\in X$, we can choose a
$G_x$-invariant $n$-form $\omega_x$ near $x$. The canonical bundle
$K_{X/G}$ admits then as a local generator the form $\omega_{X/G}$
such that $q^*\omega_{X/G}=\omega_x$, where $q$ is the quotient
map. Next assume that a crepant resolution $\pi:Y\rightarrow X/G$
exists. This means exactly that in a neighbourhood of
$\pi^{-1}(y)$,  $y:=q(x)$, the $n$-form $\pi^*\omega_{X/G}$,
defined on the smooth locus of $Y$, extends to a holomorphic
$n$-form $\omega_Y$ which  generates the canonical bundle of $Y$.
We now claim that the graph $\Gamma$ of the meromorphic map
$$q':X\dashrightarrow Y$$
is a $K$-isocorrespondence. Indeed, choose as before
$x,\,\omega_x$. Then the equality
$$q^*\omega_{X/G}=\omega_x,$$
and the fact that $\pi^*\omega_{X/G}=\omega_Y$ on the smooth locus
of $Y$ show that on the smooth part of $\Gamma$, we have
$$pr_1^*\omega_x=pr_2^*\omega_Y.$$
Since $\omega_Y$ generates $K_Y$, this shows immediately that
$\Gamma$ is a $K$-isocorrespondence.

 The simplest example of such a situation is the case of an
 involution $i$ acting with isolated fixed points on a surface $X$.
 Then the involution acts on the blow-up $\tilde X$ of $X$ at the
 fixed points, and the lifted involution
 $\tilde i$ fixes the exceptional curves pointwise.
 Then the quotient map
 $$\tilde X\rightarrow \tilde X/\tilde i$$
 ramifies simply along the exceptional curves. Furthermore
 the ramification divisor of the blowing-down map
 $\tau:\tilde X\rightarrow X$ is also the union of the exceptional
 curves with multiplicity $1$.
\end{ex}

\section{Calabi-Yau varieties and $K$-correspondences}
We consider projective $n$-dimensional complex manifolds with
trivial canonical bundle (Calabi-Yau manifolds). We shall denote
by $\omega_X$ a generator for $H^0(X,K_X)$. It can be normalized
up to a complex coefficient of modulus $1$ in such a way that
$$\Omega_X=(-1)^{\frac{n(n-1)}{2}}i^n\omega_X\wedge\overline\omega_X$$
has integral $1$ on $X$. This $\Omega_X$ is a canonical volume
form on $X$. We want to show the existence of
self-$K$-isocorrespondences for a large set of Calabi-Yau
manifolds, which have furthermore the property, like isogenies of
abelian varieties, of multiplying the canonical volume form by a
real number $>1$. The Calabi-Yau varieties for which we are able
to prove this fall into three classes. Consider the following
properties :
\begin{enumerate}
\item\label{A} $X$ is covered by abelian varieties.
\item\label{B}  There exists a rationally connected variety $Y$,
such that some embedding $j:X\hookrightarrow Y$ realizes $X$ as a
member of the linear system $\mid-K_Y\mid$ on $Y$.
\item\label{C} $X$ is the Fano variety (assumed to be smooth of the right dimension)
 of linear subspaces
${\mathbb P}^r\subset M$ of a complete intersection $M\subset
{\mathbb P}^N$ of type $(d_1,\ldots,d_k)$, with the exception of
the case $(d_1,\ldots,d_k)=(2,\dots,2)$. (Here the numbers $r,\,N,
d_i$ are chosen in such a way that $K_X$ is trivial. For fixed $r,
d_i$'s, this happens in exactly one dimension $N$ (see below).)
\end{enumerate}
We shall show later that the varieties in class \ref{C} do not in
general  satisfy property  stated in \ref{B}.

Our result is the following :
\begin{theo} \label{thm1} Assume $X$ satisfies \ref{A}, \ref{B} or
is generic satisfying \ref{C}. Then there exists a
self-$K$-isocorrespondence
$$\Sigma\subset X\times X$$
which satisfies the property that
\begin{eqnarray}\label{dilvol}f^*\Omega_X=\lambda
g^*\Omega_X,\,\lambda>1.
\end{eqnarray}
\end{theo}
Here as always, $f$ and $g$ denote the two projections to $X$, on
a desingularisation of $\Sigma$. We can rephrase formula
(\ref{dilvol}) as follows : since $\Sigma$ is a
self-$K$-isocorrespondence, there is a non zero coefficient $\mu$
 such that
 $$f^*\omega_X=\mu g^*\omega_X.$$
 Indeed, because
 $\omega_X$ nowhere vanishes, these two $n$-forms have the same zero divisor on
 $\widetilde\Sigma$, which is equal to $R_f=R_g$. So the statement concerning the
 volume form is simply the statement that we can
 find such a $\Sigma$ with $\lambda=\mid\mu\mid^2\not=1$.
 Notice that this is also equivalent to the fact that the degrees
 of $f$ and $g$ are not equal. Indeed, we have a formula
 $$f^*\Omega_X=\lambda
g^*\Omega_X$$ and $\lambda$ is then computed by integrating both
sides over $\widetilde\Sigma$, which gives
\begin{eqnarray}\label{label}deg\,f=\lambda deg\,g.
\end{eqnarray}

In case \ref{A}, the construction of $\Sigma$ is straightforward.
Namely, let
$$
\begin{matrix} &P&\stackrel{\varphi}{\rightarrow}& X\\
&h \downarrow&&\\
&B&&
\end{matrix}
$$
be a covering of $X$ by abelian varieties.  We may assume that
there is a rational section  of $h:P\rightarrow B$. Hence the
smooth fibers of $h$ have a zero, which allows to define
multiplication by any integer $m\in{\mathbb Z}$. Now choose two
integers $m$ and $m'$ and define
$$\Sigma=\{(\phi(mx),\phi(m'x)),\,x\in P\}.$$
(To be more rigorous, take the closure of the set above defined
for $x\in P^0$, the open set of $P$ where $h$ is of maximal rank.)
 It is easy to see that $\Sigma $ has dimension  $n$. The fact
 that
 it is a self-$K$-isocorrespondence follows as above, using the fact that $K_X$ is trivial,
  from the following
 formula (\ref{abvol}), where $a$ is the dimension of the abelian varieties $P_b$ :
 \begin{eqnarray}
 \label{abvol}{\frac{1}{m^a}pr_1^*\omega_X}_{\mid
 \Sigma}={\frac{1}{{m'}^a}pr_2^*\omega_X}_{\mid \Sigma}
 \end{eqnarray}
 as $n$-forms on the smooth locus of $\Sigma$.
The formula (\ref{abvol}) also  shows that the coefficient
$\lambda$ introduced above is equal to $(\frac{m}{m'})^{2a}$,
hence can be made different from $1$.

To prove formula (\ref{abvol}), we note that $\Sigma$ is the image
under $(\phi,\phi)$ of $\Sigma'\subset P\times_BP$,
$$\Sigma'=\{(mx,m'x), \,x\in P\}.$$
 Next the restriction
$\Sigma'_b$ of $\Sigma'$ to $P_b\times P_b$ is the graph
$$\Sigma'_b=\{(mx,m'x), \,x\in P_b\}.$$
It is obvious that it satisfies
$${\frac{1}{m^a}pr_1^*\omega_{P_b}}_{\mid \Sigma'_b}=
{\frac{1}{{m'}^a}pr_2^*\omega_{P_b}}_{\mid \Sigma'_b},$$ where
$\omega_{P_b}$ is a holomorphic $a$-form on $P_b$. Now, since
$\Sigma'\subset P\times_BP$, the two projections $pr_1$, $pr_2$
from $\Sigma'$ to $P$ induce
\begin{eqnarray}\label{ven}pr_1^*:R^0h_*K_{P/B}\rightarrow
K_{\Sigma'/B},\,pr_2^*:R^0h_*K_{P/B}\rightarrow K_{\Sigma'/B},
\end{eqnarray} and the maps
\begin{eqnarray}\label{ven1}pr_1^*:R^0h_*K_{P}\rightarrow
K_{\Sigma'/B},\,pr_2^*:R^0h_*K_{P}\rightarrow K_{\Sigma'}
\end{eqnarray} are simply the above tensorized with the identity of
$K_B$. Since we just noticed that the maps $pr_i^*$ in (\ref{ven})
satisfy the relation $\frac{1}{m^a}pr_1^*=\frac{1}{{m'}^a}pr_2^*$,
it follows that the same relation holds for the maps $pr_i^*$ of
(\ref{ven1}). Taking global sections, it follows that
$$\frac{1}{m^a}pr_1^*\omega_{P}=\frac{1}{{m'}^a}pr_2^*\omega_{P}$$
for any holomorphic $n$-form $\omega_{P}$ on $P$ and in particular
for $\phi^*\omega_X$. \cqfd
 {\bf Proof of
Theorem \ref{thm1} in case \ref{B}.} The construction is the
following : recall that we have an embedding
$$j:X\hookrightarrow Y.$$
Now choose a rational curve $C\subset Y$ with ample normal bundle,
which exists since $Y$ is rationally connected (\cite{komomi}). We
may assume furthermore that the intersection of $C$ with $X$, as a
divisor on $C$, is of the form
$$mx+m'y+z,$$
where $x\not=y$ and $z$ is a reduced zero-cycle on $C$ disjoint
from $x$ and $y$. Here $m$ and $m'$ are two distinct fixed
integers, on which the choice of $C$ will of course depend.   Now
choose a hypersurface $W\subset X$ supporting $z$. We will then
define, for an adequate choice of $W$, the correspondence $\Sigma$
as the closure of the image  in $X\times X$ via the map $(F,G)$
defined below, of the following set
$$\Sigma'=\{(x',y',C'),\,C'\cdot X=mx'+m'y'+z',\,z'\subset W\},$$
where in this definition, $C'$ has to be a deformation of $C$. The
map
$$(F,G):\Sigma'\rightarrow X\times X$$
is defined by
$$(F,G)((x',y',C'))=(x',y').$$
More precisely, we will consider below the (unique) component of
$\Sigma'$ passing through $(x,y,C)$. We first show that
$dim\,\Sigma'=n$ for a generic choice of $W$.  This is an easy
dimension count : the Hilbert scheme of $C$ is smooth at $C$ and
has dimension
$$h^0(C,N_{C/Y})=-K_Y\cdot C+n-2.$$
Next we impose the conditions that the intersection of $C'$ with
$X$ is finite (this is open) and of the form $mx'+m'y+z'$, with
$z'\subset W$. This imposes at most $(m-1)+(m'-1)+deg\,z'$
conditions to the deformations of $C'$.  Furthermore, one sees
easily that for an adequate choice of $W$, these conditions are
infinitesimally  independent at our initial point $(x,y,C)$. Hence
it follows that $\Sigma'$ is smooth at $(x,y,C)$, of dimension
$$dim\,\Sigma'=-K_Y\cdot C+n-2-((m-1)+(m'-1)+deg\,z')$$
$$=-K_Y\cdot C+n-2-(X\cdot C-2)$$
and this is equal to $n$ because $X\in \mid -K_Y\mid$.

An easy infinitesimal computation involving the fact that the
normal bundle of $C$ can be choosen arbitrarily ample will show
that $\Sigma$ is also of dimension $n$, or more precisely has a
component of dimension $n$.

It remains now to show that $\Sigma$ gives a
self-$K$-isocorrespondence satisfying furthermore the condition
\begin{eqnarray}\label{ello}m^2f^*\Omega_X={m'}^2g^*\Omega_X.
\end{eqnarray}
 (Choosing then $m\lneq m'$, will give a coefficient
$\lambda=\frac{{m'}^2}{m^2}\gneq1$. The formula (\ref{ello}) and
the fact that $\Sigma$ is a $K$-isocorrespondence will follow from
the following fact :
\begin{lemm} We have
\begin{eqnarray}\label{midi}mF^*\omega_X+m'G^*\omega_X=0
\end{eqnarray} on $\Sigma'$ for any holomorphic $n$-form on $X$.
\end{lemm}
Indeed, since $(\Sigma, (pr_1,pr_2))$ is the Stein factorization
of $(\Sigma',(F,G))$, the formula will be true as well for
$(\Sigma',F,G) $ replaced with $(\Sigma,pr_1,pr_2)$ or better by
 a desingularization $(\widetilde\Sigma,f,g) $.  Now, the canonical bundle of $X$ being
trivial, the divisor of $f^*\omega_X$ (resp. $g^*\omega_X$) is
equal to $R_f$ (resp. $R_g$), so that formula
\begin{eqnarray}\label{midi10}mf^*\omega_X+m'g^*\omega_X=0
\end{eqnarray}
implies that $\Sigma$ is a $K$-isocorrespondence. \cqfd
 {\bf Proof of the lemma.} We have three $0$-correspondences
 between $\Sigma'$ and $X$. The first one is
 $\Gamma_C\subset \Sigma'\times X$, which has for fiber
 over $\sigma=(x,y,C)\in \Sigma'$ the $0$-dimensional
 subscheme $C\cap X$ of $X$. If
 $$
 \mathcal{C}\subset \Sigma'\times Y$$
 is the universal subscheme, corresponding to the map from
 $\Sigma'$ to the Hilbert scheme of curves in $Y$,
 then
 $$\Gamma_C=\mathcal{C}\cap(\Sigma'\times X).$$
 The second one is
 $\Gamma_{x,y}$, whose fiber over $\sigma=(x,y,C)\in\Sigma'$ is
 the $0$-cycle $mx+m'y$. This correspondence is nothing but
 the sum
 $m\Gamma_F+m'\Gamma_G$ of the graphs of $F$ and $G$.
The third one, which we denote by $\Gamma_z$ has for fiber over
$\sigma=(x,y,C)\in\Sigma'$ the residual cycle $z=C\cdot X-mx-m'y$.
Hence we obviously have the relation
$$\Gamma_C=\Gamma_z+\Gamma_{x,y}$$
 as n-cycles in
$\Sigma'\times X$. It follows from this that for $\omega_X\in
H^0(X,K_X)$, the Mumford pull-backs $\Gamma_C^*\omega_X$,
$\Gamma_{x,y}^*\omega_X$ and $\Gamma_z^*\omega_X$ satisfy the
relation
\begin{eqnarray}\label{jm}
\Gamma_C^*\omega_X=\Gamma_z^*\omega_X+\Gamma_{x,y}^*\omega_X.
\end{eqnarray}
Since $\Gamma_{x,y}=m\Gamma_F+m'\Gamma_G$, we have
$$\Gamma_{x,y}^*\omega_X=mF^*\omega_X+m'G^*\omega_X,$$
and hence (\ref{jm}) gives
$$mF^*\omega_X+m'G^*\omega_X=\Gamma_C^*\omega_X-\Gamma_z^*\omega_X.$$
To prove (\ref{midi}), it suffices now to prove that
$\Gamma_C^*\omega_X$ and $\Gamma_z^*\omega_X$ vanish.

For the second one, this is quite easy. Indeed, by definition of
$\Sigma'$, the cycle $\Gamma_z$ is supported on $\Sigma'\times W$.
On the other hand the $n$-form $\omega_X$ vanishes on $W$, because
$dim\,W<n$. So $\Gamma_z^*\omega_X=0$.

As for the second one, we already noticed the fact that
\begin{eqnarray}\label{ilest}\Gamma_C=(Id,j)^*\mathcal{C},
\end{eqnarray}
 where we see $\mathcal{C}$ as a codimension $n$-cycle in
$\Sigma'\times Y$. This last cycle induces a cohomological
correspondence
$$[\mathcal{C}]^*:H^1(Y,\Omega_Y^{n+1})\rightarrow
H^0(\Sigma',K_{\Sigma'}).$$ Formula (\ref{ilest}) then shows
immediately that
$$\Gamma_C^*\omega_X=[\mathcal{C}]^*(j_*\omega_X).$$
So to conclude the proof that $\Gamma_C^*\omega_X=0$, it suffices
to see that $j_*\omega_X=0$ in $H^1(Y,\Omega_Y^{n+1})$. But this
last space is in fact $0$, because it is Serre dual to
$H^n(Y,\mathcal{O}_Y)$ and  $Y$ is rationally connected. \cqfd

{\bf Proof of Theorem \ref{thm1} in case \ref{C}.} The
construction in this case is as follows. We assume the complete
intersection $M\subset \mathbb{P}^N$ is a generic complete
intersection of multidegree $d_1\leqslant\ldots\leqslant d_k$, so
that its Fano variety $X$ of $r$-planes is smooth of the right
dimension. Let
$$G=Grass(r+1,N+1).$$
Then the canonical bundle $K_G$ is equal to $-(N+1)L$, where
$L=det\,\mathcal{E}$ is the Pl\"ucker line bundle, $\mathcal{E}$
is the dual of the tautological subbundle. Now $X\subset G$ is
defined as the $0$-set of the section
$(\tilde\sigma_1,\ldots,\tilde\sigma_k)$ of the vector bundle
$S^{d_1}\mathcal{E}\oplus\ldots\oplus S^{d_k}\mathcal{E}$
corresponding to the section $(\sigma_1,\ldots,\sigma_k)$ of
$\mathcal{O}_{\mathbb{P}^N}(d_1)\oplus\ldots\oplus\mathcal{O}_{\mathbb{P}^N}(d_k)$
defining $M$.

It follows from adjunction that the canonical bundle of $X$ is
given by the formula
$$K_X=-(N+1)L_{\mid X}+\sum_idet\,S^{d_i}\mathcal{E}.$$
We use the following lemma :
\begin{lemm}
Let $E$ be a vector bundle of rank $k$. Then for any integer $l$,
we have
$$det \,S^lE\cong (det \,E)^{\otimes \alpha},$$
where $\alpha=h^0(\mathbb{P}^k,\mathcal{O}_{\mathbb{P}^k}(l-1))$.
\end{lemm}
We conclude from this and the fact that
$$rk\,\mathcal{E}=r+1,\,det\,\mathcal{E}=L,$$ that
the triviality of the canonical bundle of $X$ is equivalent to the
equality
\begin{eqnarray}\label{parti}
N+1=\sum_ih^0(\mathbb{P}^{r+1},\mathcal{O}_{\mathbb{P}^{r+1}}(d_i-1)).
\end{eqnarray}
Since $d_k\geq3$, we can choose two integers $m<m'$ such that
$m+m'=d_k$. Let $Z$ be the complete intersection defined by
$(\sigma_1,\ldots,\sigma_{k-1})$. We consider now
$$\widetilde\Sigma=
\{(P_1,P_2,P),\,P\cong \mathbb{P}^{r+1}\subset Z,$$
$$P_1,\,P_2\subset  M,\,P\cap M\supseteqq mP_1+m'P_2\}.$$
(Note that $M\subset Z$ is defined by one equation of degree $d_k$
so that we have then  either $P\cap M=mP_1+m'P_2$ or $P\subset
X$.)

We define $f$ and $g$ from $\widetilde\Sigma$ to $X$ by
$$f(P_1,P_2,P)=P_1\in X,\,g(P_1,P_2,P)=P_2\in X.$$
\begin{lemm}
We have $dim\,\widetilde\Sigma=n=dim\,X$ and $f,\,g$ are
dominating.
\end{lemm}

The variety $W$ parametrizing the $\mathbb{P}^{r+1}$'s contained
in $Z$ is the $0$-set of the natural section $$
(\tilde\sigma'_1,\ldots,\tilde\sigma'_{d-1})$$ of the bundle
$S^{d_1}\mathcal{E}'\oplus\ldots\oplus S^{d_{k-1}}\mathcal{E}'$ on
the Grassmannian  $Grass(r+2,N+1)$. By genericity of $Z$, it is
smooth of dimension \begin{eqnarray} \label{dent}
(r+2)(N-r-1)-\sum_{i\leq k-1}rk\,S^{d_i}\mathcal{E}',
\end{eqnarray} where
$$rk\,S^{d_i}\mathcal{E}'=h^0(\mathbb{P}^{r+1},\mathcal{O}_{\mathbb{P}^{r+1}}(d_i)).$$

On $W$ there is a $\mathbb{P}^{r+1}\times\mathbb{P}^{r+1}$-bundle,
whose fiber over $P\in W$ parametrizes  pairs of hyperplanes in
$P$. Let us call it $W'$. Then we have
$$\widetilde\Sigma\subset W',$$
and $\widetilde\Sigma$ is defined by the condition that
$(P_1,P_2,P)\in \widetilde\Sigma$ if and only if the restriction
${\sigma_{k}}_{\mid P}$ is proportional to $\tau_1^m\tau_2^{m'}$,
where $\tau_i$ are linear equations defining $P_i$ in $P$.
 In other words, $\widetilde\Sigma\subset W'$ is the zero locus of the section
 of the vector bundle
 $\pi^*S^{d_k} \mathcal{E}'/H$, where $\pi$ is the projection from
 $W'$ to $W$ and $H$ is the line subbunble with fiber
$<\tau_1^m\tau_2^{m'}>$ at $(P_1,P_2,P)$. It follows that
\begin{eqnarray} \label{denti}
dim\,\widetilde\Sigma\geqslant dim\,W+2(r+1)-rk\,S^{d_k}
\mathcal{E}'+1, \end{eqnarray} and since our equations are
generic, a standard argument shows that we have in fact equality
and that $\widetilde\Sigma$ is smooth.
 Combining (\ref{dent}) and (\ref{denti}), we get
 $$dim\,\widetilde\Sigma=(r+2)(N-r-1)-\sum_{i}rk\,S^{d_i}\mathcal{E}'+2(r+1)+1.$$
 Next we note that, since $X$ is the zero locus of a transverse section
 of the vector bundle $\oplus_iS^{d_i}\mathcal{E}$ on
 $Grass(r+1,N+1)$,  we have
 $$dim\,X=(r+1)(N-r)-\sum_i rk\,S^{d_i}\mathcal{E}.$$
 Noting finally that if $E$ is of rank $r+1$, then
 $rk\,S^kE=h^0( \mathbb{P}^r, \mathcal{O}_{\mathbb{P}^r}(k))$,
 we get
 $$dim\,\widetilde\Sigma-dim\,X$$
 $$=
 -(\sum_i h^0( \mathbb{P}^{r+1}, \mathcal{O}_{\mathbb{P}^{r+1}}(d_i))
 -h^0( \mathbb{P}^r,
 \mathcal{O}_{\mathbb{P}^r}(d_i)))+(r+2)(N-r-1)-(r+1)(N-r)+2r+3$$
 $$=-\sum_i h^0( \mathbb{P}^{r+1}, \mathcal{O}_{\mathbb{P}^{r+1}}(d_i-1)
 +(r+2)(N-r-1)-(r+1)(N-r)+2r+3.$$
 Using equality (\ref{parti}), this gives us
 $dim\,\widetilde\Sigma-dim\,X=0$.

 To conclude the proof of the lemma, we have to show that
 the maps $f,\,g$ are dominating. We do it for $f$: the fiber of
 $f$ over $P_1\in X$ is the zero locus of a section $s$
 of a vector bundle
 over the variety $W'_{P_1}$ parametrizing $ \mathbb{P}^{r+1}$'s
 containing $P_1$, together with a hyperplane $P_2$ in them. Precisely, this vector bundle
 has for fiber
 $$\oplus_{i<k}
 H^0(P, \mathcal{O}_P(d_i-1))\oplus H^0(P,
 \mathcal{O}_P(d_k-1))/<\tau_1^{m-1}\tau_2^{m'}>$$
 at $(P,P_2)$. The section $s$ takes the value
 $$(({\sigma_i}_{\mid
 P}/\tau_1)_{i<k},{\sigma_{k}}_{\mid
 P}/\tau_1\,{\rm mod}\, \tau_1^{m-1}\tau_2^{m'}$$
 at $(P,P_2)$, where $\tau_i$ is a defining
 equation
 for $P_i\subset P$. We use here the fact that $P_1\subset M$, so that
 ${\sigma_i}_{\mid P}$ vanishes along $P_1$.
 The
 vector bundle has the same rank as the variety $W'_{P_1}$, as
 shows the previous computation. To show that $f$ is dominating,
 it suffices to show that this vector bundle has a non zero top
 Chern class, which is not hard.
 \cqfd

 To conclude the proof of the theorem in case \ref{C}, it remains
 to prove the following
 \begin{lemm} \label{petitlemme}The two ramification divisors $R_f$ and $R_g$ are
 equal.
 \end{lemm}
 Indeed, the lemma shows that $\widetilde\Sigma$ provides a
 self-$K$-isocorrespondence of $X$. Next,
 we have explained after the statement of the theorem that
 for a self-$K$-isocorrespondence of a $K$-trivial variety,
 the fact that it multiplies the volume by a real coefficient
 different from $1$ as in formula (\ref{dilvol}) is equivalent to
 the fact that the degrees of the maps $f$ and $g$ are different
 (cf (\ref{label})). Now the degree of $f$ and $g$ are the top
 Chern classes of the vector bundles described above. From this it is easy
 to show that for
 $m>m'$ we have $deg\,f<deg\,g$.

\cqfd

{\bf Proof of lemma \ref{petitlemme}.} We observe first that the
set $K$ of $(P_1,P_2)\in \Sigma$ such that the linear space
generated by $P_1$ and $P_2$ is a $ \mathbb{P}^{r+1}$ contained in
$M$ is of dimension $<n$. Suppose that we show that $R_f=R_g$ away
from $(f,g)^{-1}(K)$ : then $R_f-R_g$ is a divisor which is
rationally equivalent to $0$ (since both $R_f$ and $R_g$ are
members of the linear system $K_{\widetilde\Sigma}$), and
supported on $(f,g)^{-1}(K)$ which is contracted by $(f,g)$. But
it is well known that the components of a contractible divisor are
rationally independent. Hence this suffices to imply that
$R_f=R_g$. Next, we show by a dimension count (recall that our
parameters are generic) and the description given above of the
fibers of $f$ and $g$ that, away from $(f,g)^{-1}(K)$, the
ramification of $f$ and $g$ is simple, i.e. the ramification
divisor is reduced. In conclusion, it suffices to show that we
have the set theoretic equality $R_f=R_g$ away from
$(f,g)^{-1}(K)$. Next we note that the set of $(P_1,P_2,P)\in
\widetilde\Sigma$ such that $P_1=P_2$ is of codimension greater
than $1$ in $\widetilde\Sigma$. Hence it suffices to show the set
theoretic equality $R_f=R_g$ away from $(f,g)^{-1}(K)$ and at
points $(P_1,P_2,P)$ where $P_1\not=P_2$.

We do it by an explicit computation : let $P_1\in X$, and let
$(P_1,P_2,P),\,P_1\not=P_2,\,P\not\subset M$ be a point of
$\widetilde\Sigma$ where $f$ ramifies. So there is a first order
deformation $P_\epsilon$ of $P$ in $Z$, fixing $P_1$, and  such
that ${\sigma_k}_{\mid P_\epsilon}$ remains to first order of the
form $\tau_{1,\epsilon}^m\tau_{2,\epsilon}^{m'}$, where
$\tau_{1,\epsilon}$ is a defining equation of $P_1$ in
$P_\epsilon$. Since the first order deformation $P_\epsilon$ fixes
$P_1$, it is contained in a $ \mathbb{P}^{r+2}$ that we shall
denote by $P'$. Let us choose coordinates $X_0,\ldots, X_{r+2}$ on
$P'$ so that $P$ is defined by $X_{r+2}=0$, $P_1$ is defined by
$X_{r+2}=X_{r+1}=0$ and $P_2$ is defined by $X_{r+2}=X_{r}=0$.

The deformation $P_\epsilon$ is then given by the equation $$
X_{r+2}=\varepsilon X_{r+1}.$$ We have by assumption : $$
{\sigma_k}_{\mid P'}=X_{r+1}^mX_r^{m'}+X_{r+2}G.$$ It follows that
in the coordinates $X_0,\ldots,X_{r+1}$ for $P_\epsilon$, we have
to first order in $\epsilon$ :
$${\sigma_k}_{\mid P_\epsilon}=X_{r+1}^mX_r^{m'}+\epsilon X_{r+1}G',$$
where $G'$ is the restriction of $G$ to $P$.
 Since $\tau_{1,\epsilon}$ is proportional to $X_{r+1}$, the condition that
${\sigma_k}_{\mid P_\epsilon}$ remains to first order of the form
$\tau_{1,\epsilon}^m\tau_{2,\epsilon}^{m'}$ is then clearly
$$G'=X_{r+1}^{m-1}X_r^{m'-1}A$$
for some linear form $A$ on $P$. Hence our condition is that
\begin{eqnarray}\label{def}{\sigma_k}_{\mid
P'}=X_{r+1}^mX_r^{m'}+X_{r+1}^{m-1}X_r^{m'-1}X_{r+2}A+X_{r+2}^2H\\
=X_{r+1}^{m-1}X_r^{m'-1}(X_{r+1}X_r+X_{r+2}A)\,{\rm mod}
X_{r+2}^2.
\end{eqnarray}
Furthermore, we note that the fact that $P$ has a deformation in
$P'$ which remains contained in $Z$ can be written as the fact
that ${\sigma_i}_{\mid P'},\,i<k$ vanish at order $2$ along $P$,
hence it does not depend on $P_1$. Now the equation (\ref{def}) is
symmetric in $P_1$ and $P_2$. This shows that $g$ ramifies as well
at $(P_1,P_2,P)$, and concludes the proof of lemma
\ref{petitlemme}.

 \cqfd
As a remark of independent interest, we show now that many
varieties in class \ref{C}  do not fall in class \ref{B}.
\begin{prop} Let $X$ be the Fano variety of $ \mathbb{P}^r$'s in a
 generic complete intersection  $M$  of type $(d_1,\ldots, d_k)$
 in $ \mathbb{P}^N$. Assume
 that $X$ has trivial canonical bundle and that all $d_i$ are $>2$.
Then $X$ is not a member of $\mid-K_Y\mid$ in a rationally
connected variety $Y$.
\end{prop}
{\bf Proof.} We  note to begin with that, for such an $X$,
$Pic\,X= \mathbb{Z}L$, where $L$ is the Pl\"ucker line bundle.
Assuming the statement of the proposition is false, there would
then be two possibilities :
\begin{enumerate}
\item\label{i} $(-K_Y)_{\mid X}=\alpha L,\,\alpha\leq 0$.
\item\label{ii}$(-K_Y)_{\mid X}=\alpha L,\,\alpha> 0$.
\end{enumerate}
In case \ref{i}, $Y$ is not Fano, but since it is rationally
connected, it admits a Mori contraction $$\phi:Y\rightarrow Y'.$$
Since curves contracted by $\phi$ have negative intersection with
$K_Y$, we see that $\phi_{\mid X}$ has to be finite. It follows
that the positive dimensional fibers of $\phi$ (which all meet $X$
since $X\in \mid-K_Y\mid$), are $1$-dimensional. Hence they are
rational curves $C$ such that their space of deformation in $Y$ is
at most $n$-dimensional. There are two possibilities for generic
$C$:

- $h^0(N_{C/Y})=n$, that is $C\cdot X=2$, or

- $h^0(N_{C/Y})=n-1$, that is $C\cdot X=1$.

In the first case, $\phi_{\mid X}$ is of degree $2$, and since it
is finite, $X$ has an involution over $Y'$. We show in lemma
\ref{findepreuve} that $X$ has no non-trivial automorphism, so
this case is ruled out.

In the second case, $\phi$ is the contraction of a divisor. Then
because $X\in \mid-K_Y\mid$, one can show that the variety $Y'$ is
smooth and that $X\hookrightarrow Y'$ is realized via $\phi$ as a
member of $\mid-K_{Y'}\mid$. Induction on the Picard number of $Y$
also rules out this case.

It remains next to consider case \ref{ii}. In this case, since $
\mathcal{O}_X(X)$ is very ample, and $H^1(Y, \mathcal{O}_Y)=0$,
$X$ moves in $Y$. On the other hand, it does not deform in $Y$,
because  one can compute that $H^1(X,T_X\otimes
L^{-\alpha})=0,\,\forall\alpha>0$. If  $X$ has no non-trivial
automorphism, it follows by taking a pencil in $Y$ that $Y$ is
birationally equivalent to $X\times \mathbb{P}^1$, hence is not
rationally connected. So we get a contradiction using the
following lemma : \begin{lemm}\label{findepreuve} $X$ being as the
statement of the proposition, $X$ has no non-trivial automorphism.
\end{lemm} \cqfd {\bf Proof of lemma \ref{findepreuve}.} An
automorphism $\iota$ of $X$ induces a automorphism $\tilde\iota$
of the projective space generated by $X$ in the Pl\"ucker
embedding, because $Pic\,X={\mathbb Z}L$, where $L$ is the
Pl\"ucker line bundle. Now, one shows that, under our assumptions,
this projective space is the Pl\"ucker space, that is there are no
linear equations on the Grassmannian which vanish on $X$. Next one
shows that the quadratic equations vanishing on $X$ also vanish on
the Grassmannian. It follows that the automorphism $\tilde \iota$
of the Pl\"ucker projective space induces an automorphism of the
Grassmannian. Hence it is induced by an automorphism $\iota'$  of
the projective space $ \mathbb{P}^n$, which has to fix the variety
$M$. One concludes then using the fact that $M$ has no non trivial
automorphism.

\cqfd
\section{Intrinsic pseudo-volume forms and a problem of \\Kobayashi}
The Kobayashi-Eisenman pseudo-volume form $\Psi_X$ on a complex
manifold $X$ is defined as follows : for $x\in X,\,u\in T_{X,x}$
put
$$\Psi_X(u)=\frac{1}{\lambda},$$
where $\lambda=Max_{\phi:D^n\mapsto
X,\,\phi(0)=x}\,\{\mid\mu\mid,\,\phi_*(\frac{\partial}{\partial
z_1}\wedge\ldots\wedge\frac{\partial}{\partial z_n})=\mu u\}$.

Here $D$ is the unit disk in $ \mathbb{C}$.
\begin{rema} A similar definition can be made using the ball
instead of the polydisk , cf \cite{demailly}. The resulting
pseudo-volume forms so obtained are equivalent, and all the
results that follow will be true as well with this definition of
$\Psi_X$.
\end{rema}
 Denoting $\kappa_n$ the hyperbolic form on the polydisk :
\begin{eqnarray}\label{kappa}\kappa_n=i^n\Pi_{1}^{n}\frac{1}{(1-\mid
z_j\mid^2)^2}dz_j\wedge\overline{dz}_j, \end{eqnarray} we see
immediately, using the fact that $ \kappa_n$ coincides with the
standard volume form at $0$, is invariant under the automorphisms
of $D_n$, and that the later act transitively on $D^n$, that we
can express $\Psi_X$ as follows :
\begin{eqnarray}\label{psix}
\Psi_{X,x}=inf_{\phi:D^n\rightarrow
X,\,\phi(b)=x}\{(\phi_b^{-1})^*\kappa_n\}.
\end{eqnarray}
 Here we consider only
the holomorphic maps $\phi:D^n\rightarrow X$  which are unramified
at $b,\, \phi(b)=x$, and $\phi_b$ is then defined as the local
inverse of $\phi$ near $b$. It is obvious from either definition
that $\Psi_X$ satisfies the decreasing volume property with
respect to holomorphic maps :

{\it For any holomorphic map $\phi:X\rightarrow Y$ between
$n$-dimensional complex manifolds, we have
$$\phi^*\Psi_Y\leq \Psi_X.$$}

Also, the following theorem is a consequence of  Ahlfors-Schwarz
lemma (cf \cite{demailly}) :
\begin{theo} \label{thahl}If $X$ is isomorphic to $D_n$ (resp. to the quotient
of $D_n$ by a group acting freely and properly discontinuously, eg
$X$ is a product of curves), then $\Psi_X=\kappa_n$, (resp. to the
hyperbolic volume form on the quotient induced by $ \kappa_n$).
\end{theo}

There is also a meromorphic version $\widetilde{\Psi}_X$ which has
the advantage of being  invariant under birational maps : namely
put \begin{eqnarray}\label{tildepsix}
\widetilde\Psi_{X,x}=inf_{\phi:D^n\dashrightarrow
X,\,\phi(b)=x}\{(\phi_b^{-1})^*\kappa_n\}.
\end{eqnarray}
Here we consider  the meromorphic maps $\phi:D^n\dashrightarrow X$
which are defined at $b$ and unramified at $b, \,\phi(b)=x$, and
$\phi_b$ is then defined as the local inverse of $\phi$ near $b$.

The following result is proved in \cite{grigre},
\cite{kobaochiai2} :
\begin{theo}\label{GrGr} If $X$ is a projective complex manifold
which  is of general type, then $\widetilde{\Psi}_X$ is non
degenerate outside a proper closed algebraic subset of $X$.
\end{theo}
(The result is proved in \cite{grigre} for $\Psi_X$, but it is not
hard to see that it applies as well to $\widetilde{\Psi}_X$.)
Kobayashi \cite{kobayashi} conjectures the converse to this
statement :
\begin{conj} \label{conjko}If $X$ is a projective complex manifold
which  is  not of general type, then $\widetilde{\Psi}_X=0$ on a
dense Zariski open set of $X$.
\end{conj}
\begin{rema} A priori, $\widetilde{\Psi}_X$ is only
uppersemicontinuous, hence the equality $\widetilde{\Psi}_X=0$ on
a dense Zariski open set of $X$ does not imply that
$\widetilde{\Psi}_X=0$ everywhere.
\end{rema}
This conjecture is known in dimension $\leq2$, \cite{grigre}. In
dimension $2$, it uses the classification of surfaces, and the
fact that $K3$-surfaces are covered by elliptic curves. The proof
shows more generally that $ \Psi_X=0$ on a dense Zariski open set,
for a variety which is covered by abelian varieties.

We start this section with the definition of modified versions
$\Phi_X,\,\Phi_{X,an}$ of $\Psi_X$, together with their
meromorphic counterparts $\widetilde \Phi_X,\,\widetilde
\phi_{X,an}$.
\begin{Defi} We put
$$\Phi_{X,x}=inf_{\Sigma\subset X\times X,K-iso,\,\sigma\in\widetilde\Sigma,\,g(\sigma)=x}
(f^*\Psi_{X})_\sigma.$$
\end{Defi}
Here $\Sigma$ runs through the self-$K$-isocorrespondences of $X$,
and we denote as usual
$$\begin{matrix}&\widetilde \Sigma&\stackrel{g}{\rightarrow}&X\\
&f\downarrow&&\\
&X&&
\end{matrix}
$$
a desingularization. We use then the fact that $\widetilde\Sigma$
induces a canonical isomorphism
$$f^*K_X\cong g^*K_X$$
to see that $(f^*\Psi_X)_\sigma$ gives a pseudo-volume element for
$X$ at $x,\,g(\sigma)=x$.

The definition of $\Phi_{X,an}$ is similar:
\begin{Defi}
We put
$$\Phi_{X,an,x}=inf_{\Sigma\subset Y\times X,K-corresp,\,\sigma\in\widetilde\Sigma,\,g(\sigma)=x}
(f^*\Psi_{Y})_\sigma.$$
\end{Defi}
Here $\Sigma$ runs through the  set of all  $K$-correspondences
from $Y$ to $X$, but we ask that $\Sigma$ is unramified at
$\sigma$, namely that near $\sigma\in \widetilde\Sigma$, we have
the equality $R_f=R_g$.   Then exactly as above,
$(f^*\Psi_Y)_\sigma$ gives a pseudo-volume element for $X$ at
$x,\,g(\sigma)=x$.

The reason for the notation $\Phi_{X,an}$ ({\it an} for analytic)
is the fact that $\Phi_{X,an}$ is much more analytic, namely it
can also be computed as
\begin{eqnarray}\label{eqmar}\Phi_{X,an,x}=inf_{\Sigma\subset D^n\times
X,K-corresp,\,\sigma\in\widetilde\Sigma,\,g(\sigma)=x}
(f^*\kappa_n)_\sigma. \end{eqnarray}

\begin{rema}
There are other intermediate possible definitions for a modified
version of $\Psi_X$ using $K$-correspondences. For example, we
could restrict in the definition of $\Phi_{X,an}$ to the proper
$K$-correspondences, i.e. those for which $g$ is also proper. In
the definition of $\Phi_X$, we could consider all
$K$-isocorrespondences from $Y$ to $X$, instead of the
self-$K$-isocorrespondences from $X$ to $X$. We restricted to the
two extremal cases, which seem to be the most interesting, because
on one side $\Phi_X$ is of course the closest to $\Psi_X$, while
on the other side $\Phi_{X,an}$ satisfies the following version of
the decreasing volume property, as follows immediately from its
definition.
\end{rema}
\begin{lemm} If $\Sigma\subset Y\times X$ is a $K$-correspondence,
then with the notations used before for the desingularization :
$$g^*\Phi_{X,an}\leq f^*\Phi_{Y,an}.$$

\end{lemm}
\cqfd
 Note also that from the definition of $\Phi_X$ we get the
following :
\begin{lemm} \label{utileplustard}If $\Sigma\subset X\times X$ is a
self-$K$-isocorrespondence, we have with the same notations :
$$f^*\Phi_X=g^*\Phi_X.$$
\end{lemm}
\cqfd Finally, we define the meromorphic versions $\widetilde
\Phi_X,\,\widetilde \Phi_{X,an}$ by the formula :
$$\widetilde
\Phi_{X,x}=inf_{\phi:X\dashrightarrow
Y,\,\phi(x)=y}\,\{\phi^*\Phi_{Y,y}\},$$
$$\widetilde
\Phi_{X,an,x}=inf_{\phi:X\dashrightarrow
Y,\,\phi(x)=y}\,\{\phi^*\Phi_{Y,an,y}\}.$$ In both formulas,  we
consider only the birational maps $\phi:X\dashrightarrow Y$ which
are defined at $x$ and such that $\phi^{-1}$ is defined at
$y=\phi(x)$.

Of course we have, for any birational map $\phi:X\dashrightarrow
Y$, the equalities
$$\phi^*\widetilde
\Phi_Y=\Phi_X,$$
$$\phi^*\widetilde
\Phi_{Y,an}=\Phi_{X,an},$$ which are satisfied on the open set $U$
of $X$ where $\phi$ is defined and is a local isomorphism. In
particular, if $U\hookrightarrow X$ is the inclusion of a Zariski
open set, we have \begin{eqnarray} \label{mar22}
{{\widetilde \Phi}}_{X\mid U}=\Phi_U,\, {{\widetilde \Phi}}_{X,an\mid U}=\Phi_{U,an}.
\end{eqnarray}

Our main result towards the comparison of $\Phi_X,\,\Phi_{X,an}$
and $\Psi_X$ is the following  :
\begin{theo} \label{thsec3}If $X$ is the polydisk $D^n$, or any quotient of the polydisk
by a free properly discontinuous action of a group on $D^n$, eg
$X$ is a product of curves, then
$$
\Phi_{X,an}=\Psi_X.$$ Since
$$\Phi_{X,an}\leq\Phi_X\leq\Psi_X,$$
it follows that $\Phi_X=\Psi_X$ too.
\end{theo}

{\bf Proof.} We do it for $D^n$, the general case follows exactly
in the same way, using the fact that $\Psi_X$ in this case is the
hyperbolic volume form, which satisfies the K\"ahler-Einstein
equation (\ref{KE}). The proof is very similar to the proof that
$\Psi_{D^n}=\kappa_n$ (theorem \ref{thahl}), namely it uses the
Ahlfors-Schwarz lemma. We want however to explain carefully why it
works as well in the context of $K$-correspondences.

By formula \ref{eqmar}, what we have to prove is the following :

{\it If $$ \begin{matrix}&\widetilde \Sigma&\stackrel{g}{\rightarrow}& D^n\\
&f\downarrow&&\\
&D^n&&\end{matrix} $$ is the desingularization of a
$K$-correspondence from $D^n$ to itself, then
\begin{eqnarray} \label{mar222}g^*\kappa_n\leq
f^*\kappa_n.
\end{eqnarray}
} Let
$$ \begin{matrix}&\widetilde \Sigma_\epsilon&\stackrel{g_\epsilon}{\rightarrow}& D^n\\
&f_\epsilon\downarrow&&\\
&D^n&&\end{matrix} $$ be the restriction of the $K$-correspondence
$\Sigma$ to the polydisk of radius $1-\epsilon$. In other words we
intersect $\Sigma $ with $D^n_{1-\epsilon}\times D^n$ and we
identify $D^n_{1-\epsilon}$ with $D^n$ via the homothety of
coefficient $\frac{1}{1-\epsilon}$. It suffices to show that
\begin{eqnarray} \label{mar2222}g_\epsilon^*\kappa_n\leq
f_\epsilon^*\kappa_n.
\end{eqnarray}
Next because $ \Sigma$ is a $K$-correspondence, the ratio
$$\psi_\epsilon:=\frac{g_\epsilon^*\kappa_n}{f_\epsilon^*\kappa_n}$$
is a non negative $C^\infty$-function (which is even real
analytic) on $\widetilde \Sigma_\epsilon$. Furthermore, as we have
restricted to $D^n_{1-\epsilon}$,  the numerator stays bounded
near the boundary while the denominator tends to $\infty$
generically on the boundary of $D^n$, so we have
$$\lim_{f(x)\rightarrow\partial D^n}\psi_\epsilon(x)=0.$$
 It follows then from the properness of the map $f_\epsilon$ that
$\psi_\epsilon$ has a maximum on $\widetilde \Sigma_\epsilon$. Let
$\psi_\epsilon(x)$ be maximum. Formula (\ref{mar2222}) is
equivalent to
$$\psi_\epsilon(x)\leq1.$$
Assume the contrary and let $c:=\psi_\epsilon(x)>1$. Choose
$\alpha$ generic, $1<\alpha<c$. Let
$$\widetilde
\Sigma_{\epsilon,\alpha}:=\{y\in \widetilde \Sigma_\epsilon,\,
\psi_\epsilon(x)\geq\alpha\}.$$
 Then since $\alpha$ is generic,
and $\psi_\epsilon$ tends to $0$ near $\partial \widetilde
\Sigma_{\epsilon}$,
 $\widetilde
\Sigma_{\epsilon,\alpha}$ is compact and has a smooth boundary.
Now let $\chi=i^n\Pi_{j=1}^{j=n}\frac{1}{(1-\mid z_j\mid^2)^2}$.
By definition of $\kappa_n$ (cf (\ref{kappa})), we have
$$\psi_\epsilon=\frac{g_\epsilon^*\chi}{f_\epsilon^*\chi}\mid G\mid^2,$$
where $G$ is  holomorphic. Furthermore, we have the
K\"ahler-Einstein equation \begin{eqnarray} \label{KE}
(\frac{i}{2}\partial\overline\partial log\,\chi)^n=n!\kappa_n.
\end{eqnarray}
Denoting by $\omega=\frac{i}{2}\partial\overline\partial
log\,\chi$, we have
\begin{eqnarray}\label{avantdermar}\frac{i}{2}\partial\overline\partial
log\,\psi_\epsilon=g_\epsilon^*\omega-f^*_\epsilon\omega,\,
\omega^n=\kappa_n.\end{eqnarray} Now, in $\widetilde
\Sigma_{\epsilon,\alpha}$, we have $\psi_\epsilon>1$, which
implies that
\begin{eqnarray}\label{dermar}f^*_\epsilon\kappa_n\leq g_\epsilon^*\kappa_n,
\end{eqnarray}
with strict inequality away from the ramification divisor
$R_f$. Let
$$\theta:=g_\epsilon^*\omega^{n-1}+g_\epsilon^*\omega^{n-2}f^*_\epsilon\omega+\ldots
+f^*_\epsilon\omega^{n-1}.$$ This is a semipositive
$(n-1,n-1)$-form, which is positive away from $R_f$. Furthermore
formulae (\ref{avantdermar})  and (\ref{dermar}) say that
\begin{eqnarray}\label{finmar}
(\frac{i}{2}\partial\overline\partial
log\,\psi_\epsilon)\theta\geq0
\end{eqnarray}
in $\widetilde \Sigma_{\epsilon,\alpha}$ with strict inequality
away from the ramification divisor $R_f$. Of course, if we knew
that $x\not\in R_f$ then we would conclude that the hypothesis
that  $log\, \psi_\epsilon$ has  a maximum at $x$ is absurd,
because its Hessian should then be seminegative at $x$,
contradicting the strict inequality in (\ref{finmar}). In general,
one can apply the following (standard) argument  : choose a number
$\alpha'$, such that $\alpha<\alpha'<log\,c$. Put
$$\mu^+=Sup\,(0,log\,\psi_\epsilon-\alpha').$$
Then $\mu^+$ is non negative, vanishes  identically near the
boundary of $\widetilde \Sigma_{\epsilon,\alpha}$, and is positive
at $x$. Now consider
$$\int_{\widetilde
\Sigma_{\epsilon,\alpha}}\mu^+(\frac{i}{2}\partial\overline\partial
log\,\psi_\epsilon)\theta.$$ This is strictly positive. On the
other hand, integration by parts, using the fact that the
derivatives of  $\mu^+$ are integrable, gives :
$$\int_{\widetilde
\Sigma_{\epsilon,\alpha}}\mu^+(\frac{i}{2}\partial\overline\partial
log\,\psi_\epsilon)\theta=-\int_{\widetilde
\Sigma_{\epsilon,\alpha}}\frac{i}{2}(\partial
\mu^+\wedge\overline\partial log\,\psi_\epsilon)\theta.$$ But
since $\mu^+ =log\,\psi_\epsilon-\alpha'$ when it is non zero, the
integral on the right is equal to $$-\int_{\widetilde
\Sigma_{\epsilon,\alpha'}}\frac{i}{2}(\partial
log\,\psi_\epsilon\wedge\overline\partial
log\,\psi_\epsilon)\theta,$$ where
$$\widetilde
\Sigma_{\epsilon,\alpha'}=\{y\in \widetilde
\Sigma_{\epsilon,\alpha},\,log\,\psi_\epsilon(y)\geq\alpha'\}.$$
But this last integral is obviously negative, which is a
contradiction.

\cqfd Next we have  the following strenghtening of theorem
\ref{GrGr}.
\begin{theo} If $X$ is a projective complex manifold which is of
general type, we have $\Phi_{X,an}>0$ (and in particular
$\phi_X>0$) away from a proper closed algebraic subset of $X$.
\end{theo}

{\bf Proof.} We just sketch the argument, since it is a
combination of the construction in \cite{grigre},
\cite{kobaochiai} and of the arguments given above in the specific
case of $K$-correspondences.

Since $X$ is of general type, there exists an inclusion of sheaves
$$ L\subset K_X^{\otimes\alpha}$$
for sufficiently large $\alpha$, where $L$ is an ample line bundle
on $X$. Then, if $h_L$ is a hermitian metric on $L$ such that the
associated Chern form
$$\omega_{L,h}=\frac{1}{2i\pi}\partial\overline\partial h_L$$
is a K\"ahler form, we can see
$\mu:=\frac{1}{h_L^{\frac{1}{\alpha}}}$ as a pseudovolume form on
$X$, vanishing along a divisor, which satisfies, in local
coordinates where $\mu=i^n\chi dz_1\wedge d\overline
z_1\wedge\ldots\wedge dz_n\wedge d\overline z_n$, the equation
\begin{eqnarray}\label{omegalh1} i\partial\overline\partial log\,
\chi=\frac{1}{\alpha}\omega_{L,h}.
\end{eqnarray} Now after a rescaling, we may
assume that
\begin{eqnarray}\label{omegalh}(\frac{1}{2\alpha}\omega_{L,h})^n\geq n!\mu.
\end{eqnarray}
So the theorem is a consequence of the following proposition,
which is proved exactly as theorem \ref{thsec3} :
\begin{prop} Assume $X$ is equipped with a pseudo-volume form
$\mu$ satisfying equations (\ref{omegalh1}) and (\ref{omegalh}).
Then for any $K$-correspondence $\Sigma\subset D^n\times X$, we
have
$$g^*\mu_n\leq f^*\kappa_n.$$
\end{prop}
 \cqfd
\begin{rema} One can show similarly that the same result holds for
$\widetilde\Phi_{X,an}$.
\end{rema}
The two theorems above obviously lead to the following
\begin{conj} \label{conj}$\Phi_{X}$ is equivalent to $\Psi_X$.
This means that there exists a non zero constant $\alpha$
depending on $X$ such that
$$ \alpha\Psi_X\leq\Phi_X\leq\Psi_X.$$
\end{conj}

We conclude this section with the proof of the following theorems,
which prove a number of special cases of Kobayashi's conjecture
\ref{conjko} for our pseudovolume forms
$\Phi_X,\,\widetilde\Phi_X$.
\begin{theo}\label{thm322} Assume $X$ is a $K$-trivial projective variety
which is as in the statement of theorem \ref{thm1}, that is
satisfies \ref{A}, \ref{B} or is generic satisfying \ref{C}. Then
$\Phi_X=0$.
\end{theo}
\begin{theo} \label{thm422} Assume $X$ is birational to $X'$, and there exists a
projective morphism $\phi:X\rightarrow B$ such that
$dim\,B<dim\,X$ and the generic fiber $X'_b$ is a $K$-trivial
variety as in the previous theorem. Then $\widetilde\Phi_X=0$ on a
dense Zariski open set of $X$.
\end{theo}
{\bf Proof of theorem \ref{thm322}.} By theorem \ref{thm1}, there
exists a self-$K$-isocorrespondence $\Sigma\subset X\times X$ such
that, with the notation $$\begin{matrix}&\widetilde \Sigma&\stackrel{g}{\rightarrow}&X\\
&f\downarrow&&\\
&X&&
\end{matrix}
$$ for a desingularization of $\Sigma$, we have
$$f^*\Omega_X=\lambda g^*\Omega_X,$$ for
some $\lambda>1$. By lemma \ref{utileplustard}, we know that
$$f^*\Phi_X= g^*\Phi_X.$$
Writing $\Phi_X=\chi\Omega_X$ and combining these two equalities
gives
$$f^*\chi=\lambda g^*\chi.$$
But the function $\chi$ is uppersemicontinuous and bounded, hence
it has a maximum. Let $x$ be a point where $\chi(x)$ is maximum.
Let $\sigma\in \widetilde\Sigma$ be such that $g(\sigma)=x$. Then
for $y=f(x)$, we get $\chi(y)=\lambda \chi(x)$. Since $\chi(x)$ is
maximum, we also have $\chi(y)\leq \chi(x)$, which implies
$\chi(x)=0$ because $\lambda>1$. So $\chi=0$.

\cqfd

{\bf Proof of theorem \ref{thm422}.} By the birational invariance
of $\widetilde\Phi_X$, it suffices to show that
${\widetilde\Phi_{X'}}$ on a dense Zariski open set $X''$ of $X'$,
or equivalently that, for some dense open $X''$, one has
$\widetilde\Phi_{X''}=0$.

But the construction of the self-$K$-isocorrespondence given in
the proof of theorem \ref{thm1} can be made in families at least
over a Zariski open set $B''$ of $B$. Letting
$X''=\phi^{-1}(B'')$, we get a relative self-$K$-isocorrespondence
$$\Sigma\subset X''\times_{B''}X''$$
satisfying the property that as relative pseudovolume forms on
$\widetilde\Sigma$ over $B''$, we have
\begin{eqnarray}\label{etcest}f^*\Omega_{X''/B}=\lambda g^*\Omega_{X''/B},\,\lambda>1.
\end{eqnarray}
(Indeed, note that the coefficient $\lambda$ is constant in
families, by the formula (\ref{label}).) Now, since $\Sigma\subset
X''\times_{B''}X''$, we have $\phi\circ f=\phi\circ g=:\pi$ and
for $\Omega_B$ a volume form on $B$, we have
$$f^*(\Omega_{X''/B}\otimes\phi^*\Omega_B)=(f^*\Omega_{X''/B})\otimes\pi^*\Omega_B,$$
and similarly for $g$. Hence (\ref{etcest}) gives
\begin{eqnarray}\label{etcestencore}f^*(\Omega_{X''/B}\otimes\phi^*\Omega_B)=\lambda
g^*(\Omega_{X''/B}\otimes\phi^*\Omega_B).
\end{eqnarray} Denoting by $\mu$ the
volume form $\Omega_{X''/B}\otimes\phi^*\Omega_B$ on $X''$, we
have a relation $\Phi_{X''}=\chi\mu$ for some function $\chi$, and
formula (\ref{etcestencore}), together with the relation
$$f^*\Phi_{X''}= g^*\Phi_{X''}$$
show that
$$f^*\chi=\lambda g^*\chi.$$
One concludes then that $\chi=0$, hence $\Phi_{X''}=0$, using the
fact that $\Sigma\subset X''\times_BX''$ and the properness of
$\phi:X''\rightarrow B''$.

\cqfd

\section{Concluding remarks and questions}
\subsection{Fano varieties of $r$-planes in a hypersurface}
Our first question concerns the Chow-theoretic interpretation of
our construction of a self-$K$-correspondence in case \ref{C},
that is when $X$ is the variety of $r$-planes in a hypersurface of
degree $d$ (or more generally a complete intersection). Unlike
case \ref{B}, we did not deduce formula
 \begin{eqnarray}\label{nov1}f^*\omega_X=\mu g^*\omega_X
 \end{eqnarray}
 from a relation between $0$-cycles, of the shape
 \begin{eqnarray}\label{nov2}\forall \sigma\in\Sigma,\,\alpha f(\sigma)+\beta
 g(\sigma)\equiv z,
 \end{eqnarray}
 where $z$ is supported on a proper algebraic subset of $X$, and
 $\alpha,\,\beta$ are fixed integers depending on
 the integers $m,\,m'$. Of course, by Mumford's theorem \cite{mumford}, (\ref{nov2})
 implies (\ref{nov1}), with $\mu=\frac{-\beta}{\alpha}$.
 Bloch-Beilinson's conjectures predict also that conversely
(\ref{nov1})
 implies (\ref{nov2}). So our first question is : how to prove a
 formula like (\ref{nov2}), for $\Sigma$ constructed as in
 the proof of theorem \ref{thm1}, case \ref{C}?
 Let us do it in the case where $M$ is the  cubic fourfold, $m=2,\,m'=1$ and
 $r=1$. In this case $X$ is $4$-dimensional and is hyperK\"ahler
 (cf \cite{beauvilledonagi}).

 Recall that $\Sigma$ parametrizes the pairs $(L_1,L_2)$ of lines
 in
 $M$ such that there exists a plane $P\subset \mathbb{P}^5$, with
 $$P\cap M=2L_1+L_2.$$
 For each line $L\subset M$, let us denote by $l$ the corresponding
 point in
 $X$. For a generic $l\in X$ there is an incidence surface in $X$
 (cf \cite{voisin})
 $$S_l:=\{l'\in X,\,L\cap L'\not=\emptyset\}.$$
 Note that if
 $$\begin{matrix}
 &P&\stackrel{q}{\rightarrow}& M\\
 &p\downarrow&&\\
 &X&&
 \end{matrix}$$
 is the incidence correspondence, we have
 $$S_l=p_*q^*L\, {\rm in}\, CH^2(X).$$
 It follows that, denoting by $h$ the class of a plane section of
 $M$, we have for any $(l_1,l_2)\in \Sigma$ the relation
\begin{eqnarray}\label{nov3}
2S_{l_1}+S_{l_2}=p_*q^*h\,{\rm in}\, CH^2(X).
 \end{eqnarray}
 Now, it is not hard to prove the following
\begin{lemm}\label{lenov}
 There exists an integer $\alpha\not=0$ and a proper algebraic
 subset $Z\subset X$ such that for any $l\in X$ the following
 relation holds in $CH(X)$ :
\begin{eqnarray}\label{nov4}
S_l^2=\alpha l+z,
 \end{eqnarray}
 where $z$ is a $0$-cycle supported on $Z$.
 \end{lemm}
We now combine formulas (\ref{nov4}) and (\ref{nov3}) to get for
any $(l_1,l_2)\in \Sigma$ the relations
$$4S_{l_1}^2=S_{l_2}^2+z',$$
$$4\alpha{l_1}=\alpha{l_2}+z'+z'',$$
where $z'$ and $z''$ are supported on a fixed algebraic subset of
$X$. This gives us the formula (\ref{nov2}) in this case. This
also shows that $\mu=\frac{1}{4}$ hence
$\lambda=\frac{1}{16}=\frac{deg\,f}{deg\,g}$ in this case.
\subsection{Some examples satisfying the Kobayashi conjecture}

 In a different direction we observe that our construction in
 case \ref{C} provides for $d=3$ a true rational map
 $\phi:X\dashrightarrow X$. Here we consider
 the Fano variety of $r$-planes in a hypersurface $M$
 of degree $3$ in $\mathbb{P}^n$, with the relation
 \begin{eqnarray}\label{reparti}n+1=h^0( \mathbb{P}^{r+1},
 \mathcal{O}_{\mathbb{P}^{r+1}}(2))
 \end{eqnarray}
 which implies that $K_X$ is trivial (cf (\ref{parti})). Now let
 as in  section 2
 $$\Sigma=\{(P_1,P_2)\in X\times X,\,\exists
 P\subset\mathbb{P}^n,\, P\cap M=2P_1+P_2\}.$$
 Here $P$ has to be a $r+1$-plane and  we in fact have to consider
 the Zariski closure of the set above.

 We have the following
 \begin{lemm} The first projection $pr_1:\Sigma\rightarrow X$
 is of degree $1$. Hence $\Sigma$ is the graph of a rational map
 $\phi$.
 \end{lemm}
 {\bf Proof.}
Let $P_1$ be generic in $X$. Consider
$$P:=\bigcap_{x\in P_1}T_{X,x}.$$
Here $T_{X,x}$ is the projective hyperplane tangent to $X$ at $x$.
Then $P$ has dimension $n-h^0(( \mathbb{P}^{r},
 \mathcal{O}_{\mathbb{P}^{r}}(2))$ because the Gauss map of
 $M$ is given by polynomials of degree $2$.
 But we have by (\ref{reparti})
 $$n-h^0(( \mathbb{P}^{r},
 \mathcal{O}_{\mathbb{P}^{r}}(2))=
-1+h^0( \mathbb{P}^{r+1},
 \mathcal{O}_{\mathbb{P}^{r+1}}(2))-h^0(( \mathbb{P}^{r},
 \mathcal{O}_{\mathbb{P}^{r}}(2))$$
 $$=-1+h^0( \mathbb{P}^{r+1},
 \mathcal{O}_{\mathbb{P}^{r+1}}(1))=r+1.$$
 Hence $P$ is a $\mathbb{P}^{r+1}$ everywhere tangent to $M$ along
 $P_1$. Hence we have
 $$P\cap M=2P_1+P_2,$$
 for some $P_2$ which must be the only point in the fiber
 of $\Sigma$ over $P_1$.

 \cqfd
\begin{coro} \label{coro}For such $X$, we have
$$\widetilde \Psi_X=0$$
on a Zariski open set of $X$. In other words, Kobayashi's
conjecture \ref{conjko} is true for $X$.
\end{coro}
{\bf Proof.} The decreasing volume property for $\widetilde
\Psi_X$ will say that
$$ \phi^*\widetilde \Psi_X \leq \widetilde \Psi_X $$
on the open set where $\phi$ is defined. On the other hand, we
have seen that $$\phi^*\Omega_X=\lambda
\Omega_X,\,\lambda=deg\,\phi,$$ where $\Omega_X$ is the canonical
volume form of $X$.
 Now we conclude as in the proof of Theorem \ref{thm322}, using the fact that
 $deg\,\phi>1$, (for example $deg\,\phi=16$ in the case of the cubic fourfold).  \cqfd
\begin{rema} The existence of the self-map $\phi$ of  degree $>1$,
hence multiplying the volume form by a coefficient $>1$, suggests
that not only the Kobayashi pseudovolume form of $X$ vanishes but
also the Kobayashi pseudodistance of $X$ vanishes, as conjectured
in \cite{campana}, \cite{kobayashi}. This would follow, as the
following argument shows, from a dynamical study of the map $\phi$
but we have not been able to do it. In fact, what is easily seen
is the fact that the Kobayashi pseudodistance $d_K$ of $X$ as
above is $0$ if, for general $y\in X$, the orbit
$\{\phi^k(y),\,k\in{\mathbb Z}\}$ is dense in $X$. Indeed, one
sees easily that $\phi$ has one fixed point $x$. Next consider the
function $\chi(y)=d_K(x,y)$ on $X$. By the decreasing distance
property, we have
$$d_K(x,\phi(y))\leq d_K(x,y).$$
So it follows that we have the inequality of pseudo-volume forms :
$$\phi^*(\chi\cdot\Omega_X)\leq \chi\cdot\phi^*\Omega_X.$$
Now we have
$$\phi^*\Omega_X=deg\,\phi\cdot \Omega_X.$$
So
$$\phi^*(\chi\cdot\Omega_X)\leq deg\,\phi \cdot \chi\cdot\Omega_X.$$
But the integrals of both sides over $X$ are equal. Hence we
conclude that
$$f^*\chi=\chi$$
almost everywhere on $X$. So we have $d_K(x,\phi(y))=d_K(x,y)$ for
almost all $y$. So if the $\phi^k(y)$ are dense in $X$, (for $k$
negative or positive), hence arbitrary close to $x$, we find that
$d_K(x,y)=0$ for almost every $y$.
\end{rema}
\subsection{$K$-correspondences and the Kodaira dimension}
 The following two propositions relate the Kodaira dimension and
 $K$-correspondences.
 \begin{prop} \label{pronov1} Let
 $\Sigma\subset Y\times X$ be a $K$-correspondence, where $X$ and
 $Y$ are smooth and projective. Then
 $$\kappa(Y)\geq\kappa(X).$$

 \end{prop}
 {\bf Proof.}
 Let
 $$\begin{matrix}&\widetilde \Sigma&\stackrel{g}{\rightarrow}&X\\
&f\downarrow&&\\
&Y&&
\end{matrix}
$$
be a desingularization of $\Sigma$. If $\kappa(X)=-\infty$ there
is nothing to prove. If $\kappa(X)\geq 0$ there is a non zero
section of $K_X^{\otimes m}$ for some $m\geq1$. Since
$g^*K_X\subset f^*K_Y$, there is a non zero section of
$f^*K_Y^{\otimes m}$, and it follows that there is a non zero
section of $K_Y^{\otimes Nm}$, where $N$ is the degree of $f$. So
$\kappa(Y)\geq0$. So we can consider the Iitaka fibration
$$Y\dashrightarrow B$$
whose generic fiber $Y_b$ satisfies
$$\kappa({K_Y}_{\mid Y_b})=0.$$
Let $\widetilde Y_b:=f^{-1}(Y_b)$. Then
$$\kappa({f^*K_Y}_{\mid
\widetilde Y_b})=0.$$ Since $ g^*K_X\subset f^*K_Y$, it follows
that
$$\kappa({g^*K_X}_{\mid
\widetilde Y_b})=0.$$ Hence the components of $g(\widetilde Y_b)$
are contained in a fiber of the Iitaka fibration $X\dashrightarrow
B'$ of $X$. It follows that $dim\,B'\leq dim\,B$.
 \cqfd
 \begin{prop}\label{pronov2}  If $X$ is a projective variety which
 is of general type, any self-$K$-isocorrespondence
 $$\begin{matrix}&\widetilde \Sigma&\stackrel{g}{\rightarrow}&X\\
&f\downarrow&&\\
&X&&
\end{matrix}
$$ satisfies
$$deg\,f=deg\,g.$$

 \end{prop}
 {\bf Proof.} For a line bundle $L$ on a projective variety, whose
 Iitaka dimension is equal to $n=dim\,X$, define
 $$d^+(L)=Sup_m\{ \frac{deg\,\phi_m(X)}{m^n}\}$$
 where $\phi_m$ is the rational map to projective space
 given by the sections of $L^{\otimes m}$ assuming there are any.
 This is a finite positive number.
 Also, it is immediate to see that if
 $$\phi:X'\rightarrow X$$
 is a generically finite cover, we have
 $$d^+(\phi^*L)=deg\,\phi \,d^+(L).$$
 We can apply this to $K_X$ and to
 $f:\widetilde \Sigma\rightarrow X$ and $g:\widetilde \Sigma\rightarrow
 X$, since
  the Iitaka dimension
 of $K_X$ is equal to $n$, and using the fact $f^*K_X\cong g^*K_X$,
 we find that
 $$deg\,f\, d^+(K_X)=deg\,g\,d^+(K_X).$$
 Hence $deg\,f=deg\,g$.
 \cqfd
 Note that in the above propositions, we  used only the fact
 that $f^*K_X\cong g^* K_X$, which is weaker than the equality
 of the ramification divisors.
 Note also that the hypothesis in proposition \ref{pronov2} is
 necessary. Indeed we know the existence of
 self-$K$-isocorrespondences $\Sigma$ of arbitrary large degree $deg\,g/deg\,f$
 for $K$-trivial varieties $X$ (eg take for $X$ an elliptic curve). Considering a product
 $Y\times X$, and the
 self-$K$-isocorrespondences $\Delta_Y\times \Sigma$ of $Y\times
 X$, we find examples of self-$K$-isocorrespondences
 of degree $\not=1$ on varieties with any possible Kodaira
 dimension, except for the maximal one.

We note also that, in the second example studied in section 2, we
have found for the considered $K$-trivial varieties $X$,
self-$K$-isocorrespondences $\Sigma\subset X\times X$ which pass
through an arbitrary pair $(x,y)$. We ask conversely:
\begin{conj}\label{conjnov} Let $X$ be such that there exists a
self-$K$-isocorrespondence $\Sigma\subset X\times X$ which passes
through an arbitrary pair $(x,y)$, then $\kappa(X)\leq0$.
\end{conj}
In case where $X$ is a curve, we have the following :
\begin{prop}\label{pronov3} Assume $X$ is a smooth curve of genus
$>1$, then any self-$K$-isocorres-\\
pondence of $X$ is rigid.
\end{prop}
{\bf Proof.} Let
$\widetilde\Sigma\stackrel{(f,g)}{\rightarrow}X\times X$ be the
desingularization of $\Sigma$. By rigid, we mean here that there
is no deformation of the triple $(\widetilde\Sigma,f,g)$, keeping
the property that $R_f=R_g$. But, since both $f$ and $g$ ramify
along $R_f$, the torsion free part of the normal bundle
$$(f,g)^*(T_{X\times X})/T_{\widetilde\Sigma}$$
has degree $\leq deg\,f^*T_X$, which is negative. Hence it has no
non zero section.

\cqfd


\begin{thebibliography}{99}
\bibitem{batyrev} V. Batyrev, Stringy Hodge numbers of varieties
with Gorenstein canonical singularities, in {\it Integrable
systems and and algebraic geometry} (Kobe/Kyoto, 1997), 1-32,
World Sci. Publishing, River edge, NJ, 1998.
\bibitem{batyrevdais} V. Batyrev, D. Dais. Strong McKay
correspondence, string theoretic Hodge numbers and mirror
symmetry, Topology 35 (1996), 901-929.
\bibitem{beauvilledonagi} A. Beauville, R. Donagi. A. Beauville,
R. Donagi, La vari{\'e}t{\'e} des droites d'une hypersurface cubique
de dimension $4$, C. R. Acad. sci. Paris, S{\'e}rie 1, 301 (1985),
703-706.
\bibitem{campana}F. Campana. Special varieties, preprint 2002.
\bibitem{demailly} J.-P. Demailly. Algebraic criteria for Kobayashi
hyperbolic projective varieties and jet differentials, in {\it
Proceedings of Symposia in Pure Mathematics volume 62.2, 1997},
285-360.
\bibitem{demaillyetal} J.-P. Demailly, L. Lempert, B. Shiffman. Algebraic
approximations of holomorphic maps from Stein domains to
projective manifolds, Duke Math. J., vol. 76, 2, 1994, 332-363.
\bibitem{denef} J. Denef, F. Loeser. Germs of arcs on singular
algebraic varieties and motivic integration, Invent. Math. 135
(1999), 201-232.
\bibitem{grigre} M. Green, P. Griffiths. Two applications of algebraic
geometry to entire holomorphic mappings, in {\it The Chern
symposium 1979}, Hsiang, Kobayashi, Singer, Weinstein Eds,
Springer-Verlag 1980, 41-74.
\bibitem{kawamata} Y. Kawamata. $D$-equivalence and
$K$-equivalence, preprint math.AG/0205287, to appear J. Diff.
Geom.
\bibitem{kobayashi} S. Kobayashi. Intrinsic distances, measures and geometric
function theory, Bull. Amer. Math. Soc. 82 (1976), 357-416.
\bibitem{kobaochiai} S. Kobayashi, T. Ochiai. Mappings into compact complex
manifolds with negative first Chern class, J. Math. Soc. Japan 23
(1971), 137-148.
\bibitem{kobaochiai2} S. Kobayashi, T. Ochiai. Meromorphic mappings into
complex manifolds of general type, Invent. Math. 31, 1975,7-16.
\bibitem{komomi} J. Koll\'ar, S. Mori, Y. Myiaoka. Rationally
connected varieties, J. Alg. Geom. 1 (1992), 429-448.
\bibitem{looijenga} E. Looijenga. Motivic measures, S{\'e}minaire
Bourbaki, 52{\`e}me ann{\'e}e, 1999-2000, n$^0$ 874.
\bibitem{mumford} D. Mumford. Rational equivalence of zero-cycles on surfaces,
J. Math. Kyoto Univ. 9 (1968), 195-204.
\bibitem{voisin} C. Voisin. Th{\'e}or{\`e}me de Torelli pour les
cubiques de $\mathbb{P}^5$, Invent. Math. 86 (1986), 577-601.
\bibitem{wang} C-L Wang. $K$-equivalence in birational geometry,
preprint 2002.

\end{thebibliography}
\end{document}